\documentclass{article}
\usepackage{graphicx}
\usepackage{amsmath}
\usepackage{amssymb}
\usepackage{hhline}
\usepackage{a4wide,epsfig}
\newtheorem{theorem}{Theorem}

\newtheorem{theoremb}{Theorem}

\newtheorem{dfn}[theoremb]{Definition}

\newenvironment{proof}[1][Proof]{\textbf{#1: }}{\vspace{3pt}}
\newtheorem{Cor}{Corollary}
\newtheorem{Lemma}{Lemma}

\newtheorem{Rem}{Remark}

\newcommand\bg{\bar g}
\newcommand\h{h_\text{\rm top}}
\renewcommand\l{\lambda}
\newcommand\po{$\!\!\!{\text{\bf.}}$ }

\newcommand\bib[1]{\bibitem[#1]{#1}}

\newcommand\C{{\mathbb C}}
\newcommand\op[1]{\mathop{\rm #1}\nolimits}
\newcommand\Q{{\mathbb Q}}
\newcommand\R{{\mathbb R}}
\newcommand\hps{\hskip-16pt . \hskip2pt}

\newcommand{\eqdef}{\stackrel{{\rm def}}{=}}

\newcommand{\diag}{\mbox{\rm diag}}

\newcommand{\Reg}{\op{Reg}}
\newcommand{\Sing}{\mbox{\rm Sing}}
\newcommand{\Supp}{\mbox{\rm Supp}}

\newcommand{\weg}[1]{}

\begin{document}

   \title{Strictly non-proportional geodesically  \\
   equivalent metrics have $h_\text{top}(g)=0$}
 \author{Boris S. Kruglikov and Vladimir S. Matveev}
 \date{}
 \maketitle

\section{\hps Definition and main results}\label{S1}

 \begin{dfn}\po
Two  ($C^\infty$-smooth) Riemannian metrics $g$ and $\bg$ on a
 manifold $M^n$ are said to be {\bf geodesically equivalent}
if their geodesics coincide as unparameterized curves. They are
{\bf strictly non-proportional at $x\in M^n$},  if  the polynomial
$\det(g_{|x}-t\bg_{|x} )$ has only simple roots. 
 \end{dfn}

The question of  whether two different metrics can have the same
geodesics is natural and, therefore, classical. The first examples
are due to E. Beltrami \cite{B}, a local descriptions of
geodesically equivalent metrics was understood  by  U. Dini
\cite{Di} and T. Levi-Civita \cite{LC}. We will recall
Levi-Civita's Theorem in Section~\ref{preliminaries1}. For more
historical  details, see the surveys \cite{Mi,Am}, or/and the
introductions to the papers \cite{M1,M4}.

The main result of our paper is (for definition and properties of
$\h$ we refer to \cite{Bo,KH,Ma}):
 \begin{theorem}\po \label{th1}
Suppose  the Riemannian metrics  $g$ and $\bg$ on a closed
connected manifold  $M^n$ are  geodesically equivalent and
strictly non-proportional at least at one point. Then the
topological entropy $\h(g)$ of the geodesic flow  of $g$ vanishes.
 \end{theorem}

The condition that the metrics are strictly non-proportional is
important: for example, the product metric on a closed product
manifold $M=M_1\times M_2$ admits a family $g_1+tg_2$ of
non-proportional metrics (but not strictly non-proportional if
$\dim M>2$) with the same geodesics. But if at least one factor
has fundamental group with positive exponential growth (for
instance if $M_1$ is  hyperbolic), then by the Dinaburg Theorem
any geodesic flow on $M$ has $\h(g)>0$.

Vanishing of the topological entropy  of a $C^\infty$-smooth flow
implies a lot of dynamical restrictions. For example, the ball
volume grows sub-exponentially with its radius (Manning's
inequality \cite{Mn}), the number of geodesic arcs joining two
generic points grows sub-exponentially with its maximal length
(Ma\~n\'e's formula \cite{Ma}) and the volume of a compact
submanifold propagated by the geodesic flow also changes
sub-exponen\-tially (Yomdin's Theorem \cite{Y}), see also
\cite{P2}.

Probably even more interesting are topological restrictions
implied by $\h(g)=0$. The subexponential growth of $\pi_1(M^n)$
(Dinaburg's Theorem \cite{D}) is not very intriguing under the
assumptions of Theorem~\ref{th1}, since it is known \cite{M3} that
in this case the fundamental group is virtually abelian. But the
restriction coming from the Gromov-Paternain Theorem \cite{G,P1} and from 
\cite{PP1}
are  new, nontrivial and interesting: Namely in the simply connected
case the manifold $M^n$ is {\bf rationally elliptic,} i.e.
$\pi_*(M^n)\otimes\Q$ is finite-dimensional. This is a very
restrictive property since by the results of \cite{FHT,Pa} a {
rationally elliptic\/} manifold $M^n$ enjoys the following
properties:

 \begin{enumerate}

\item $\dim\pi_*(M^n)\otimes\Q\le n$, $\dim H_*(M^n,\Q)\le 2^{n-1}$,  $\dim H_i(M^n,\mathbb{Q})\le
\frac{1}{2}\left(\begin{array}{c}n \\ i\end{array}\right)$ \ $(i=1,...,n-1),$
 \item The Euler characteristic $\chi(M^n)$ satisfies 
$2^n-n+1\ge \chi(M^n)\ge0$. Moreover,  
 $\chi(M^n)>0$ iff
$H_\text{odd}(M^n,\Q)=0$. 
 \end{enumerate}
A manifold $M$ with finite $\pi_1(M)$ is
called {\bf rationally hyperbolic\/},  if its universal cover  is
not rationally elliptic. Thus, as a  consequence of
Theorem~\ref{th1}, we get

 \begin{Cor}\po\label{crl1}
A rationally hyperbolic closed manifold $M^n$ does not admit two
geodesically equivalent Riemannian metrics $g$ and $\bg$  which
are strictly non-proportional at least at one point.
 \end{Cor}

Rational hyperbolithity means  nothing in dimensions less than 4,
since  all closed 4-manifolds with finite fundamental group are
rational-elliptic. Note  that  the topology of closed 2- and
3-manifolds admitting non-proportional geodesically equivalent
metrics is completely understood: In dimension 2, such manifolds
are homeomorphic to the sphere, the projective plane, the torus or
the Klein bottle \cite{MT2}. In dimension 3, such manifolds are
homeomorphic  to lens spaces or to Seifert manifolds with zero
Euler number \cite{M2}.

Starting from dimension 4, almost all simply-connected manifolds are
 rationally hyperbolic. For example, in dimension 4, up to
homeomorphism, there exist infinitely many simply-connected closed
manifolds, and only five of them are rationally elliptic: $S^4$,
$S^2\times S^2$, $\C P^2$, $\C P^2\#\C P^2$ and $\C P^2\#\overline{\C
P^2}$. It is possible to construct geodesically equivalent metrics on
$S^4$ and $S^2\times S^2$ that are strictly non-proportional at least
at one point. We conjecture here that these two are the only closed
simply-connected 4-manifolds admitting strictly non-proportional
geodesically equivalent metrics. In dimension 5, a closed
rational-elliptic manifold has rational homotopy type of $S^2\times
S^3$ or $S^5$ (there are infinitely many homotopy types for
simply-connected 5-manifolds).   By recent results of \cite{PP1} 
(see Theorem E there), a closed manifold admitting a metric with
 zero topological entropy  is  $S^5$, $S^3\times S^2$, $SU(3)/SO(3)$ or  the nontrivial $S^3$-bundle over $S^2.$
We conjecture that $S^3\times S^2$ and
$S^5$ are the only closed simply-connected connected 5-manifolds
admitting geodesically equivalent metrics which are strictly
non-proportional at least at one point.

In Section~\ref{lastsec} we  announce the  restrictions on the topology
of non-simply-connected manifolds (admitting geodesically equivalent metrics which are strictly non-proportional at least at one point) that follows from
Corollary~\ref{crl1}.

Now let us comment the proof of Theorem~\ref{th1}. The main
ingredients are Theorems~\ref{integrability}, \ref{ordered123} and
Corollary~\ref{independent}, which imply that the geodesic flow of
$g$ is Liouville-integrable.

Precisely the same integrable systems were recently actively 
studied  in  mathematical physics, in the framework  of the theory 
of separation of variables.  Depending on the school,
they are called L-systems \cite{Be},
Benenti-systems \cite{IMM} and  quasi-bi-hamiltonian systems \cite{CST}.

But Liouville integrability does not immediately imply vanishing
of the topological entropy; counterexamples can be found in
\cite{BT1,BT2,Bu1,Bu2,K,KT}. If the singularities of the integrable system
behave sufficiently good (non-degenerate in the sense of
Williamson-Vey-Eliasson-Ito \cite{E,I}, see \cite{P1},  or the Taimanov conditions
\cite{T}), or if the system has a lot of symmetries 
 (for example, as  in collective integrability \cite{BP,P1}), then $\h(g)=0$. 
But for other situations nothing is known (at least if $n>2$, see \cite{P0}), even  if  the integrals are
real-analytic or  polynomial in momenta.  

It is worth mentioning that geodesically equivalent metrics are
usually not real-analytic: Levi-Civita's Theorem from
Section~\ref{preliminaries1} shows the existence of an
infinite-dimensional space of nonanalytic $C^\infty$-per\-tur\-bations in the
class of geodesically-equivalent metrics. Also the set of singular
points of the constructed integrals for the corresponding
Hamiltonian system can be quite complicated. For instance, the
projection of the singularities in $TM^n$ to the base $M^n$ is
surjective for $n>2$ and its restriction to a singular Liouville
fiber can have image which  is  locally the product of the Cantor
set and the $(n-1)$-dimensional disk.

The logic of our proof for Theorem~\ref{th1} is as  follows:
 \begin{enumerate}
\item We show that the topological entropy is supported on the singularities,
      which we describe.
\item We show that dynamics on them can be considered as a subsystem of the
      geodesic flow
  \begin{itemize}
\item on a lower-dimensional closed submanifold
\item admitting geodesically equivalent metrics which are  strictly non-proportional at least at one point.
\end{itemize}
Therefore we can apply induction by the dimension.
\end{enumerate}

\subsection*{Acknowledgments} We thank Professors    Bangert, Butler, 
 Katok, Paternain,  Shevchishin,   Taimanov and Wilking
 for useful discussions. The second author
 thanks the  University of Tromso,
where the essential part of the results were obtained, for
hospitality, and   DFG-programm 1154 (Global Differential
Geometry) and Ministerium f\"ur Wissenschaft, Forschung und Kunst
Baden-W\"urttemberg  (Elitef\"orderprogramm Postdocs 2003) for
partial financial support.

\section{Geometry behind the geodesic equivalence}

In what follows we always assume that the manifold $M^n$ is
connected and that the Riemannian metrics $g$ and $\bar g$ on $M^n$
are geodesically equivalent and  strictly non-proportional at
least at one point.

\subsection{Integrability and Levi-Civita's Theorem} \label{preliminaries1}

A Riemannian metric $g$ determines the map $\flat_g:TM\to T^*M$
with the inverse $\sharp^g:T^*M\to TM$. Consider the (1,1)-tensor
(automorphism field) $L:TM\to TM$ given by the formula
 \begin{eqnarray}
 \label{l}
L=\left(\det(\sharp^{\bar g}\circ\flat_g)\right)^{-\frac1{n+1}}
\cdot(\sharp^{\bar g}\circ\flat_g).
 \end{eqnarray}
In local coordinates, $L^j_i=\sqrt[n+1]{(\det(\bar
g)/\det(g))}\,g_{i\alpha}\bar g^{\alpha j}$. This tensor $L$
determines the family $S_t\in C^\infty(T^*M\otimes TM)$, $t\in
\R$, of $(1,1)$-tensors
 \begin{equation}\label{st}
S_t:= \det(L-t\op{Id})\cdot(L-t\op{\rm Id})^{-1}.
 \end{equation}

 \begin{Rem}\hspace{-2mm}{\bf .}
Although $(L-t\op{Id})^{-1}$ is not defined for $t\in\op{Sp}(L) $,
the tensor $S_t$ is well-defined for every $t\in\R$. In fact, it
is the adjunct matrix of $(L-t\op{Id})$. Thus by the Laplace main
minors formula, $S_t$ is a polynomial in $t$ of degree $n-1$ with
coefficients being $(1,1)$-tensors.
 \end{Rem}

The isomorphism $\flat^g$ allows us to identify the tangent and
cotangent bundles of $M^n$. This identification allows us to
transfer the natural Poisson structure and the Hamiltonian system
$H(x,p)=\frac12p\cdot\sharp^{g}(p)$ from $T^*M^n$ to $TM^n$.

 \begin{theorem}[\cite{MT1}]\hspace{-2mm}{\bf .}
\label{integrability}
 If $g$, $\bar g$ are geodesically
 equivalent,
then, for every  $t_1,t_2\in R$, the functions
\begin{equation}\label{integral}
I_{t_i}:TM^n\to \mathbb{R}, \ \ I_{t_i}(v):= g(S_{t_i}(v),v)
\end{equation}
are commuting integrals for the geodesic flow
 of  $g$.
 \end{theorem}

Since $L$ is self-adjoint with respect to both $g$ and $\bar g$,
the spectrum $\op{Sp}(L)$ is real at every point $x\in M^n$.
Denote it by $\l_1(x)\le\dots\le\l_n(x)$. Every eigenvalue
$\lambda_i(x)$ is at least continuous functions on $M^n$, and is
smooth near the points where it is a simple eigenvalue.

 \begin{theorem}[\cite{M1}]\po\label{ordered123}
Let $(M^n, g)$ be a geodesically complete connected Riemannian
manifold. Let a Riemannian metric $\bar g$ on $M^n$ be
geodesically equivalent to $g$. Then, for every  $i\in
\{1,\dots,n-1\}$ and for all $x,y\in M^n$, the following holds:
 \begin{enumerate}
\item $\lambda_i(x)\le \lambda_{i+1}(y)$.
\item  If $\lambda_i(x)< \lambda_{i+1}(x)$,
then $\lambda_i(z)< \lambda_{i+1}(z)$
for almost every point $z\in M^n$.
\item If  $\lambda_i(x)=\lambda_j(y)$ for a certain $
j\ne i$, then there exists
$z\in M^n$ such that   $\lambda_i(z)=\lambda_j(z)$.
\end{enumerate}
 \end{theorem}

 \begin{Cor}[\cite{MT3}]\po \label{independent}
Let $(M^n, g)$ be a  connected Riemannian manifold. Suppose a
Riemannian metric $\bar g$ on $M^n$ is  geodesically equivalent to
$g$ and is strictly non-proportional to $g$ at least at one point.
Then, for every mutually-different $t_1,t_2,\dots,t_n\in\R$, the
integrals $I_{t_i}$ are functionally independent almost
everywhere, i.e. the differentials $dI_{t_i}$ are linearly
independent a.e. in $TM$.
 \end{Cor}

Let us describe the local form of the integrals $I_t$. For every
$x\in M^n$ consider  coordinates  in $T_xM^n$ such that the metric
$g$ is given by the diagonal matrix $\diag(1,1,\dots,1)$ and the
tensor  $L$ is given by the diagonal matrix
$\diag(\lambda_1,\lambda_2,\dots,\lambda_n)$. Then the tensor
(\ref{st}) reads:
 \begin{eqnarray*}
S_t& = & \det(L-t\op{Id})(L-t\op{Id})^{(-1)} \\
&=& \diag(\Pi_1(t),\Pi_2(t),\dots,\Pi_n(t)),
 \end{eqnarray*}
where the polynomials  $\Pi_i(t)$ are  given by the formula
 $$
\Pi_i(t)\eqdef \prod_{j\ne i} (\lambda_j-t) \,.
 $$
Hence, for every  $\xi=(\xi_1,\dots,\xi_n)\in T_xM^n$, the
polynomial $I_t(x,\xi)$ is given by
 \begin{equation}\label{33}
I_t=\xi_1^2\Pi_1(t)+ \xi_2^2\Pi_2(t) \dots +\xi_n^2\Pi_n(t).
 \end{equation}

For further use, let us  consider the one parameter family of
functions
 $$
I'_t\eqdef \frac{d}{dt}\bigl(I_t\bigr).
 $$
For every fixed $t\in\R$ this function is an integral of the
geodesic flow for $g$.

Let us now formulate (a weaker version of) the classical
Levi-Civita's Theorem.

 \begin{theorem}[Levi-Civita \cite{LC}]\hspace{-2mm}{\bf .}
 \label{Levi-Civita}
Consider two  Riemannian metrics on an open subset $U^n\subset
M^n$ and the tensor $L$ given by (\ref{l}). Suppose the spectrum
$\op{Sp}(L)$ is simple at every point $x\in U^n$.

Then the metrics are geodesically equivalent on $U^n$ if and only
if around each point $x\in U^n$ there exist coordinates $x_1, x_2,
\dots,x_n$ in which the metrics have the following model form:
 \begin{eqnarray}
ds_{g}^2& = & |\Pi_1(\lambda_1)|dx_1^2+
|\Pi_2(\lambda_2)|dx_2^2\ +\cdots+
 |\Pi_{n}(\lambda_n)|dx_{n}^2,\label{Canon_g}\\
ds_{\bar g}^2&=&\rho_1|\Pi_1(\lambda_1)|dx_1^2+
\rho_2|\Pi_2(\lambda_2)|dx_2^2+\cdots+
\rho_{n}|\Pi_{n}(\lambda_n)|dx_{n}^2 \label{Canon_bg},
 \end{eqnarray}
where the functions $\rho_i$ are given by
 \begin{eqnarray*}
\rho_i  & \eqdef & \frac{1}{\l_1\l_2\dots\l_n}\frac{1}{\l_i}.
 \end{eqnarray*}
and $\l_i=\l_i(x_i)$ are smooth functions of one variable.
 \end{theorem}

 \begin{dfn}\po
The above coordinates will be called {\bf Levi-Civita coordinates}
and the neighborhoods where the coordinates are defined  will be
called  {\bf Levi-Civita charts}.
 \end{dfn}

In Levi-Civita coordinates the tensor $L$ is diagonal
$\diag(\lambda_1,\dots,\lambda_n)$, so the notations in the
Levi-Civita Theorem are compatible with those in the beginning of
the section.

 \begin{Cor}[\cite{M1,BM}]\hspace{-2mm}{\bf .} \label{nijenhuis}
Suppose the Riemannian metrics $g$, $\bar g$ are geodesically
equivalent on $M$.  Then, the  Nijenhuis torsion of the tensor $L$
given by (\ref{l})  vanishes: $N_L=0$.
 \end{Cor}

If the metrics are strictly non-proportional at least at one point, Corollary~\ref{nijenhuis}
follows from the above version of Levi-Civita's theorem. In the
general case, Corollary~\ref{nijenhuis} follows from the original
version of Levi-Civita's Theorem \cite{LC} and was proven in
\cite{M1} and \cite{BM}.

 Combining formulae (\ref{Canon_g}) and
(\ref{33}), we see that in the Levi-Civita coordinates the
function $I_t$ is given by
 \begin{equation}\label{LC-I}
I_t=\sum_i|\Pi_i(\l_i(x))|\,\Pi_i(t)\,\xi_i^2
 \end{equation}
In particular, the function  $I_{\lambda_i(x)}$ as the function on
the cotangent bundle is equal to $(-1)^{i-1}p_i^2$.

\subsection{ Distributions of eigenvectors: submanifolds $M_A$ }
\label{ma}

We begin with investigation of the set of points from the
Levi-Civita charts, the union of which is the open dense set
 $$
\op{Reg}(M)=\{x\in M :\ \l_i(x)\ne\l_j(x)\text{ for }i\ne j\}.
 $$
This set can be represented as the intersection
$\op{Reg}(M)=\cap_A\op{Reg}_A(M)$ by all (proper) subsets
$A\subset \{1,2,\dots,n\}$, where we denote
 $$
\Reg_A(M)=\{x\in M : \ \forall i\in A \ \forall j\not\in A \ \
\lambda_i(x)\ne \lambda_j(x)\}.
 $$

At every point $x\in \Reg_A(M)$ denote by $D_A(x)$ the subspace of
$T_xM^n$  spanned by the eigenspaces with the eigenvalues
$\lambda_i$, where $i\in A$. Since the eigenvalues $\lambda_i$ for
$i\in A$ do not bifurcate with the eigenvalues $\lambda_j$ for
$j\not\in  A$, $D_A$ is a smooth distribution on $\Reg_A(M)$. By
Corollary~\ref{nijenhuis} it is integrable. We will denote by
$M_A(x)$ its integral submanifold containing $x\in\Reg_A(x)\subset
M^n$.

 \begin{Lemma}\po\label{M_A}
For $x\in\Reg_A(M)$ the following statements hold:
 \begin{enumerate}
\item  The restrictions of $g$ and $\bar g$ to $M_A(x)$ are geodesically
equivalent.
\item $g|_{M_A(x)}$ and $\bar g|_{M_A(x)}$ are 
 strictly non-proportional at least at one point.
\item For $i\in A$ the $i^\text{th}$ eigenvector of $L$ (corresponding to $\l_i$)
coincides with the respective eigenvector of the operator $L_A$,
constructed via (\ref{l}) for the metrics $g|_{M_A(x)}$ and $\bar
g|_{M_A(x)}$.
\item There exists a universal along $M_A(x)$ constant $c$
(calculated explicitly in the proof) such that the part of
$c\cdot\op{Sp}(L)$, corresponding to $A$, coincides with the
spectrum of the operator $L_A$, constructed by the restricted to
$M_A(x)$ metrics.
\item In particular, if  an  eigenvalue $\lambda_i$, $i\in A$  is constant, then the
corresponding   eigenvalue of the operator $L_A$, constructed for
the restrictions of  $g$ and $\bar g$ to  $M_A(x)$, is constant on
$M_A(x)$.
\end{enumerate}
 \end{Lemma}

 \begin{proof}
The distribution $D_A$ defines a foliation on $\Reg_A(M)$ and on
its open dense subset $\Reg(M)$. Then it is sufficient to prove
the first, third and the fourth statements of the lemma at the
points of this subset. By Theorems~\ref{ordered123},
\ref{Levi-Civita} in a neighborhood of every point $x\in\Reg(M)$,
there exist Levi-Civita coordinates such that the metrics  $g,\
\bg$ are given by formulas (\ref{Canon_g})-(\ref{Canon_bg}). In
these coordinates, $M_A(x)$ is  the coordinate plaque of the
coordinate collection $x_\alpha$ with $\alpha\in
A=\{\alpha_1,\dots,\alpha_m\}$. Then the restrictions of the
metrics to $M_A(x)$ are given by:
 \begin{eqnarray*}
g_{|M_A}& = & \
|\Pi_{\alpha_1}(\lambda_{\alpha_1})|dx_{\alpha_1}^2+ \
|\Pi_{\alpha_2}(\lambda_{\alpha_2})|dx_{\alpha_2}^2\ +\cdots+ \
|\Pi_{\alpha_m}(\lambda_{\alpha_m})|dx_{\alpha_m}^2,\\
{\bg_{|M_A}}&=&\rho_{\alpha_1}|\Pi_{\alpha_1}(\lambda_{\alpha_1})|dx_{\alpha_1}^2+
\rho_{\alpha_2}|\Pi_{\alpha_2}\lambda_{\alpha_2}|dx_{\alpha_2}^2+\cdots+
\rho_{\alpha_m}|\Pi_{\alpha_m}(\lambda_{\alpha_m})|dx_{\alpha_m}^2.
 \end{eqnarray*}
Since $\lambda_j$ is constant on $M_A(x)$ for every  $j\not\in A$,
every factor of $\Pi_{\alpha_i}$ of the form $\lambda_j -
\lambda_{\alpha_i}$  can be ``hidden'' in $dx_{\alpha_i}^2$.  We
see that then the first metric is already in the Levi-Civita form,
and the second metric becomes in the Levi-Civita's form after
multiplication by
 \begin{equation} \label{c}
C\eqdef \prod_{j\not\in A} \lambda_j,
 \end{equation}
which is constant on $M_A(x)$. Hence, by  Levi-Civita's Theorem,
the restrictions of the metrics to $M_A$ are geodesically
equivalent.

Direct calculations show that in local coordinates the tensor
$L_A$ is given by:
 \begin{equation}\label{L_A}
C^{1/(m+1)}\diag(\lambda_{\alpha_1},\dots,\lambda_{\alpha_m}).
 \end{equation}
The third and the fourth statements of the lemma follow.

Now let us prove the second statement. Suppose the restriction of
the metrics are not strictly non-proportional at every point of  a
certain $M_A(x)$. Then, by Theorem~\ref{ordered123}, there exist
$\alpha_1,\alpha_2\in A$  such that
$\l_{\alpha_1}\equiv\l_{\alpha_2}$ on $M_A(x)$. Consider the set
$B:=\{1,\dots,n\}\setminus A$. Take the union of all leaves $M_B$
containing at least one point of $M_A(x)$. Clearly, this union
contains an open subset of $M^n$. Since the eigenvalues
$\lambda_{\alpha_1},\ \lambda_{\alpha_2}$ are constant along
$M_B$, in view of (\ref{L_A}) and Theorem~\ref{ordered123}, at
every point of this open subset we have
$\lambda_{\alpha_1}=\lambda_{\alpha_2}$, which contradicts
Theorem~\ref{ordered123}. Lemma~\ref{M_A} is proven.
 \end{proof}

 \begin{Lemma}\po \label{compact}
Suppose the eigenvalue $\lambda_i$ is not a constant. Take a point
$y\in M^n$ such that
 $$
\max_{x\in M}\lambda_{i-1}(x)<\lambda_i(y)<\min_{x\in
M}\lambda_{i+1}(x).
 $$
 (We assume  by definition that $\min_{x\in
M}\lambda_{n+1}(x)=\infty$ and  $\max_{x\in
M}\lambda_{0}(x)=-\infty$.)

Let ${\textsf{C}(i)}:=\{1,2,\dots,n\}\setminus \{i\}$. Then,
$M_{\textsf{C}(i)}(y)$ is a closed submanifold.
 \end{Lemma}
The conditions that the eigenvalue is not constant and that
$\lambda_i$ is neither maximum nor minimum are important: one can
construct counterexamples, if one of these conditions is omitted.

 \vspace{5pt}\begin{proof}[Proof of Lemma~\ref{compact}]
Since $\max_{x\in M}\l_{i-1}(x)<\lambda_i(y)<\min_{x\in M}\l_{i+1}(x)$,
there exist  $c_{\textrm{small}},c_{\textrm{big}}\in \R$  such
that
\begin{itemize}
\item $c_{\textrm{small}}<\lambda_i(y)<c_{\textrm{big}}$,
\item at least one of  the numbers $c_{\textrm{small}},c_{\textrm{big}}$ is a
regular value of the function $\lambda_i$,
\item the other number is not a critical  value of $\lambda_i$
(i.e. is either a regular value or  is equal to  $\lambda_i$   at
no point.)
\end{itemize}

 Denote by $N$
the connected component of the set
 $$
\{ x \in M^n: \ \ c_{\textrm{small}}\le \lambda_i(x)\le
c_{\textrm{big}} \},
 $$
containing the point $y$. Then $N\subset \Reg_{\textsf{C}(i)}(M)$
is a connected manifold with boundary. Therefore,
$D_{\textsf{C}(i)}$ is a smooth distribution on $N$. Since it is
integrable by Corollary~\ref{nijenhuis}, it defines a foliation.
By Corollary~\ref{nijenhuis}, the function $\lambda_i$ is constant on the leaves of the
foliation. Then, every connected component of the boundary of $N$
is a leaf  of the foliation.

At every  $x\in M^n$, consider the vector $v_i$ satisfying
 \begin{equation} \label{v}
\left\{  \begin{array}{ccc}L(v_i)&=&\lambda_i(x)v_i \\
g(v_i,v_i)&=& |\Pi_i(\lambda_i)|. \end{array}\right.
 \end{equation}
By definition of $N$, the function $|\Pi_i(\lambda_i)|$ is nonzero
and smooth at every point of $N$. Thus $v_i$ vanishes nowhere in
$N$. Hence,   at least on the double-cover of $N$, it  is defined
globally up to a sign and is smooth. The double-cover projection
maps closed submanifolds into closed ones. Therefore,  without
loss of generality we can assume that the vector field $v_i$ is
globally defined already on $N$.

Consider the flow of the vector field $v_i$. It takes leaves to
leaves. Indeed, it is sufficient to prove this almost everywhere,
for instance in Levi-Civita charts. In Levi-Civita coordinates the
leaves of the foliation are the plaques of the coordinates
$x_\alpha$, where $\alpha\in {\textsf{C}(i)}$, and the vector field $v_i$ is
$\pm \frac{\partial }{\partial x_i}$, so the claim is trivial.

Since the leaves are $(n-1)$-dimensional and the flow of $v_i$
shuffles them, the flow acts transitively and all leaves are
homeomorphic. Every connected component of the boundary of $B$ is
compact and is a leaf, whence all leaves are compact. In
particular, $M_{\textsf{C}(i)}(y)$ is compact. Lemma~\ref{compact}
is proven.
 \end{proof}

\subsection{ Bifurcation of eigenvalues: submanifolds $\Sing_i^j$ }
\label{ma+}

The spectrum $\op{Sp}(L)$ is simple in $\Reg(M)$, i.e. almost
everywhere in $M^n$. But at certain points the multiplicity of
some $\lambda_i$ can become greater than one. Such points will be
called {\bf the bifurcation points} of $\lambda_i$. By
Theorem~\ref{ordered123} the following types of bifurcations of
the eigenvalue $\lambda_i$ are possible.

 {\bf Case 1:}
The eigenvalues $\lambda_i$ and $\lambda_{i+1}$ are not constant
and there exists  $x\in M$ such that
$\lambda_{i}(x)=\lambda_{i+1}(x)$. Denote
$\bar\l_i=\max\l_i(x)=\min\l_{i+1}(x)$. Let us consider the set
 $$
\Sing_i^1\eqdef\{ x\in M^n: \ (\lambda_{i}(x)-\bar
\lambda_i)(\lambda_{i+1}(x)-\bar \lambda_i)=0\}.
 $$
This set was studied in \cite{M1} (see Theorem 6 there). It was
shown that $\Sing_i^1$ is a connected closed totally geodesic
submanifold of codimension one. The restrictions of the metrics to
it are strictly non-proportional at least at one point. Note that
not all points of $\Sing_i^1$ are points of bifurcation of the
eigenvalues $\l_i,\l_{i+1}$.

 {\bf Case 2:}
There exists  $x\in M$ and $i\in \{2,\dots,n-1\}$ such that
$\lambda_{i-1}(x)=\lambda_{i+1}(x)$. In this case, the eigenvalue
$\lambda_i$ is constant. Let us consider the set
 $$
\Sing_i^2\eqdef\{ x\in M^n: \
(\lambda_{i-1}(x)-\lambda_i)(\lambda_{i+1}(x)-\lambda_i)=0\}.
 $$
This set was also studied in \cite{M1} (see Theorem 6 there).  It
was shown that $\Sing_i^2$ is a connected closed totally geodesic
submanifold of codimension two. The restrictions of the metrics to
it are  strictly non-proportional at least at one point. Moreover, the set of
the points $x\in \Sing_i^2$ such that $\l_{i-1}(x)=\l_{i+1}(x)$ is
nowhere dense in $\Sing_i^2$.

 {\bf Case 3a:}
The eigenvalue $\lambda_i$ is constant, there exists $x\in M$ such
that $\lambda_{i}=\lambda_{i+1}(x)$ and there exists no $y$ such
that  $\lambda_{i-1}(y)=\lambda_i$.

 {\bf Case 3b:}
The eigenvalue $\lambda_i$ is constant, there exists $x\in M$ such
that $\lambda_{i}=\lambda_{i-1}(x)$ and there exists no $y$ such
that  $\lambda_{i+1}(y)=\lambda_i$.

In Cases 3a, 3b, let us consider respectively the sets
 $$
\Sing_i^3=\{x\in M^n:\ \lambda_{i}=\lambda_{i+1}(x)\}\quad
\text{or}\quad \Sing_i^3=\{x\in M^n:\
\lambda_{i}=\lambda_{i-1}(x)\}.
 $$
 
The next lemma shows that, similar to Cases 1 and 2, $\Sing_i^3$ is a 
submanifold of codimension $2$  and  the
restrictions of the metrics to $\Sing_i^3$ are geodesically
equivalent and strictly non-proportional at least at one point. Note that,
contrast to the previous cases, the set $\Sing_i^3$ is not
necessary connected.

 \begin{Lemma}\po\label{3}
Under assumptions of Cases 3a or 3b, the set $\Sing_i^3$ is a
\begin{itemize}
\item[(1)] totally geodesic
\item[(2)] closed submanifold of codimension 2.
\item[(3)] Moreover, the restrictions of the metrics to $\Sing_i^3$   are strictly
non-proportional at least at one point.
\end{itemize}
 \end{Lemma}
Here we will proof  that   $\Sing_i^3$  is a closed submanifold of
codimension 2 such that the restrictions of the metrics to it are
strictly non-proportional at least at one point. The first
statement of  the lemma, namely that $\Sing_i^3$ is totally
geodesic, will follow immediately from Theorem~\ref{fake}, see
Remark~\ref{rk4}. Before Theorem~\ref{fake}, Lemma~\ref{3} will be
used only once, namely in the proof of Theorem~\ref{singular}.
Since the proof of Theorem~\ref{fake} does not require
Theorem~\ref{singular}, no logical loop appears.

 {\bf Proof of  statements 2,3 of Lemma~\ref{3}:}
We consider Case 3a, the other case is completely analogous. By
definition, the set $\Sing_i^3$ is closed and, therefore, compact.

Let us show that  locally $\Sing_i^3$ is a submanifold of
codimension 2. Let $A=\{i,i+1\}$. Take a point $x_0$ such that
$\lambda_i=\lambda_{i+1}(x_0)$. Then $x_0\in\Reg_A(M)$ and we can
consider the set $M_A(x_0)$. By Lemma~\ref{M_A}, the restrictions
of the metrics to $M_A(x_0)$ are geodesically equivalent and
strictly non-proportional at least at one point. Since $M_A(x_0)$
is two-dimensional, the set of points, where these restrictions
are proportional, is discrete \cite{MT2}. In view of
Lemma~\ref{M_A}, the restrictions of the metrics are proportional
at $x_0$. Then in a small neighborhood of $x_0$, there exists no
other point $x\in M_A(x_0)$ such that
$\lambda_i=\lambda_{i+1}(x)$.  Denote by $B$ the set
$\{1,2,\dots,n\}\setminus  A$. For every point $x$ of a small
neighborhood of $x_0$ in $M_A(x_0)$, consider the set $M_B(x)$. It
is a submanifold of codimension two. Since the eigenvalues
$\lambda_i,\lambda_{i+1}$ are constant along $M_B$, in a small
neighborhood of $x_0$ the set $\Sing_i^3$  coincides with
$M_B(x_0)$. Thus it is a submanifold of codimension $2$.

By  the second statement of Lemma~\ref{M_A}, the   restrictions of the metrics 
 to $\Sing_i^3$ are strictly non-proportional at least at one point.
The $2^\textrm{nd}$ and $3^\textrm{d}$ statements of  Lemma~\ref{3} are 
 proven.

 Let us note that for a fixed $i$ only one of the submanifolds
$\Sing^j_i$, $j=1,2,3$, can be non-empty.

\section{Description of singular points} \label{sec3}

Consider some mutually-different numbers $t_1,\dots,t_n\in \R$ and
the respective integrals $I_{t_1}, \dots,I_{t_n}$.
Consider the Poisson action of the the group $(\R^n,+)$ on
$TM^n$: an  element $(a_1,...,a_n)\in \R^n $ acts by time-one  shift
along the Hamiltonian vector field of the function 
$a_1I_{t_1}+...+a_nI_{t_n}$.  Since the functions are  commuting integrals, the action is well-defined, smooth, symplectic,
 preserves the integrals $I_t$ and the Hamiltonian of the geodesic flow, see \S 49 of \cite{A} for details.

A point $(x,\xi)\in TM$ is called {\bf singular} if the
differentials $dI_{t_1},\dots,dI_{t_n}$ are linearly dependent at
$(x,\xi)$. An
orbit of the action is called {\bf singular} if it has a singular
point. All points of  a singular orbit are singular and have
the same  coefficients of the linear dependence.

Although the Poisson action   depends on the choice of
constants $t_1,\dots,t_n$, the property of $(x,\xi)$ being
singular does not depend on the choice of $t_i$ as far as these
numbers are all different.

\subsection{Singular points in  Levi-Civita coordinates} \label{sing}

The next theorem describes  singular points that 
lie over a Levi-Civita chart $U^n\subset\Reg(M^n)$. Fix a point
$x\in\Reg(M^n)$ and denote by $\bar\l_1,\dots,\bar\l_n$ the
constants $\l_1(x),\dots,\l_n(x)$ respectively.

 \begin{theorem}\po \label{singular}
Let the metrics $g$ and $\bar g$ be given by formulas
(\ref{Canon_g})-(\ref{Canon_bg}) in a neighborhood $U^n\subset
M^n$. If the point $(y,\xi)=(x_1,\dots,x_m, \xi_1,\dots,\xi_m)\in
T\Reg(M^n)$ is singular, then there exists $i\in
\{1,\dots,n\}$ such that $dI_{\bar\l_i}=0$. Then
$I_{\bar\l_i}(x,\xi)=0$ and at least one of the following
statements holds:
 \begin{enumerate}
\item The  derivative $\frac{\partial \l_i(x)}{\partial x_i}$ vanishes at
$x$. 
\item The function $I'_{\bar \lambda_i}$ vanishes at $(x,\xi)$.
 \end{enumerate}
Moreover,  if $M_{\textsf{C}(i)}(y)$ is compact,  the  whole
geodesic passing through $y$  with the velocity vector $\xi$ is
contained in $M_{\textsf{C}(i)}(y)$, where $\textsf{C}(i)$ is the
same as in Lemma \ref{compact}.
 \end{theorem}
Actually, the assumption that  $M_{\textsf{C}(i)}(y)$ is compact
is not necessary: Theorem~\ref{singular} remains  true, if we
replace this condition by the condition that $y\not\in \Sing_i^1$.
Our stronger assumption makes the proof shorter.

 {\bf Proof of Theorem~\ref{singular}:} Suppose the point $(y,\xi)$ is
 singular. Then, there exist constants $(\mu_1,\dots,\mu_n)\ne (0,\dots,0)$
 such that at $(y,\xi)$ it holds:
 $$
\mu_1 dI_{\bar \lambda_1} +\dots+\mu_n dI_{\bar \lambda_n}=0.
 $$

We will show that for every $i$ such that  $\mu_i\ne 0$ the
differential $dI_{\bar \lambda_i}$ vanishes at $(y,\xi)$. For
every $j\in \{1,...,n\}$  consider the function
$I_{\lambda_j(x)}(x,\eta):=\left(I_t(x,\eta)\right)_{|t=\lambda_j(x)}$.
In a small neighborhood of $y$, the function $\lambda_j$ is
smooth. Hence the function $I_{\lambda_j(x)}$ is smooth as well.
At the point $(y,\xi)$ we have:
 $$
dI_{\lambda_j(y)}=dI_{\bar \lambda_j} + I'_{\bar \lambda_j}\cdot
d\lambda_j.
 $$
We will work on the cotangent bundle to $M^n$. As we explained in
Section \ref{preliminaries1}, the function $I_{\lambda_j(x)}$ is
equal to $(-1)^{j-1}p_j^2$ and its differential has coordinates
 $$
(\underbrace{0,\dots,0}_{n+j-1},2\cdot (-1)^{j-1}\cdot
p_j,0,\dots,0).
 $$
Since the function $\lambda_j$ depends on $x_j$ only, its
differential is
 $$
(\underbrace{0,\dots,0}_{j-1},\frac{\partial \lambda_j}{\partial
x_j},0,\dots,0).
 $$
Thus $dI_{\bar\lambda_j}$ at $(y,\xi)$ is given by
 $$
(\underbrace{0,\dots,0}_{j-1},I'_{\bar \lambda_j}\cdot
\frac{\partial \lambda_j}{\partial
x_j},\underbrace{0,\dots,0}_{n-1},2\cdot (-1)^{j-1}\cdot
p_j,0,\dots,0).
 $$

We see that the differentials $dI_{\bar \lambda_j}$ do not
combine: If $\mu_i\ne 0$, then $dI_{\bar \lambda_i}=0$. Therefore,
$p_i=0$ (i.e. $\xi_i=0$), which is equivalent to
$I_{\bar\l_i}(x,\xi)=0$, and at least one of the following holds:
$\frac{\partial \lambda_i}{\partial x_i}(x)=0$ or $I'_{\bar
\lambda_i}(x,\xi)=0$. The first part of the theorem is proven.

Now let us show that the geodesic $\gamma$ such that
$(\gamma(0),\dot\gamma(0))=(y,\xi)$ is contained in
$M_{\textsf{C}(i)}(y)$.  Since $M_{\textsf{C}(i)}(y)$ is compact,
 it is sufficient to prove that at almost every point of the
geodesic the velocity vector of the geodesic is contained in
$D_{\textsf{C}(i)}$. Since $\Sing_k^j$ are   totally geodesic
submanifolds, the geodesic $\gamma$ intersect them transversally,
and  it is sufficient to prove that the velocity vector of the
geodesic lies in $D_{\textsf{C}(i)}$ in Levi-Civita's charts.

Since $I_{\bar\lambda_i}$ is an integral and $dI_{\bar
\lambda_i}=0$ at $(y,\xi)$, we obtain that $dI_{\bar \lambda_i}$
vanishes at every point $(\gamma(t),\dot\gamma(t))$. Then, as we
explained above,  in the Levi-Civita  chart, the component $\xi_i$
equals zero, so that the velocity vector of the geodesic lies in
$D_{\textsf{C}(i)}$. Finally, the geodesic stays in
$M_{\textsf{C}(i)}$ forever. Theorem~\ref{singular} is proven.

\subsection{Removable singularities} \label{sing+}

Our next goal is to show that certain singular points  are
artificially singular: if we use a finite cover and choose the
integrals  appropriate, they become regular.

Suppose the eigenvalue $\lambda_i$ is constant. From the proof of
Theorem~\ref{singular} it follows that   for every
$x\in\Reg_{\{i\}}(M)$ and  $\xi\in D_{\textsf{C}(i)}(x)\subset
T_xM^n$  the differential $dI_{\lambda_i}$ vanishes at $(x,\xi)$.
We  will show that this singularity is {\bf removable}, in the
sense that on an appropriate finite cover we can find a linear in
velocities function $J_i$ such that
$J_i^2=(-1)^{i-1}I_{\lambda_i}$. This relation immediately implies
that $J_i$ commutes with the functions $I_t$. Since
$I_{\lambda_i}$ is an integral, $J_i$ is an integral as well.
Since it is linear in velocities,  it corresponds to a Killing
vector field. We will show that this Killing vector field is
nonzero at $x$, which automatically implies that  the differential
of this integral does not vanish at $(x,\xi)$.

 In the Levi-Civita coordinates $I_{\lambda_i}=(-1)^{i-1}p_i^2$
and  we can put $J_i=\pm p_i$. Clearly, in the Levi-Civita
coordinate system, $J_i(\eta):=g(v_i,\eta)$, where
$v_i=\pm\frac{\partial}{\partial x_i}$.

Note that  the vector field $\frac{\partial}{\partial x_i}$
satisfies conditions (\ref{v}), and  that near every regular point
every vector field satisfying   (\ref{v}) is the vector field
$\frac{\partial}{\partial x_i}$ of  a certain Levi-Civita
coordinate system.

Thus, in order to show that (at least on a finite cover) there
exists a smooth function $J_i$ such that it is linear in
velocities and such that $J_i^2=(-1)^{i-1}I_{\lambda_i}$, it is
sufficient to prove
 \begin{theorem}\po\label{fake}
Suppose $\lambda_i$ is constant. Then at least on a double cover
of $M^n$ there exists a smooth vector field $v_i$ satisfying
(\ref{v}) at every point $x\in M^n$.
 \end{theorem}
\begin{Rem}\po\label{rk4}
Conditions (\ref{v}) imply that   the zeros of $v_i$ coincide with
$\cup_{j=2,3}\Sing_i^j$. Since $v_i$ is a Killing vector field,
$\Sing_i^3$ is a totally-geodesic
 submanifold.
  \end{Rem}
 {\bf Proof of Theorem~\ref{fake}:}
First we show that at least on the double-cover there exists  a
continuous  vector field $v_i$ with the required properties. In
order to do this, it is sufficient to  prove   the following
semi-local statement:
 \begin{enumerate}
\item[\bf (S)]
Locally near every point $x$ there exist  precisely two continuous
vector fields $v_i$ satisfying (\ref{v}).
 \end{enumerate}

If $\lambda_{i-1}(x)\ne \lambda_i \ne \lambda_{i+1}(x)$, then
$y\in\Reg_{\{i\}}(M)$. Then, $\Pi_i(\lambda_i)\ne 0$. Hence,
$v_i\ne 0$ in a small  neighborhood of $x$ and  the statement {\bf
(S)} is trivial.

Let us consider $x\in \Sing_i^j$, where $j=2$ or $3$, and prove
the statement in a small disk  neighborhood $U^n\ni x$.

First of all, if a vector field $v_i$ satisfies (\ref{v}), then
the vector field $-v_i$ satisfies (\ref{v}) as well. Since
$\Sing_i$ is nowhere dense, the fields do not coincide. Therefore
we obtain at least two different required vector fields.

Next, there exist no more than two such vector fields. Indeed,
such a vector field $v_i$ must vanish along $\Sing_i^j$, since
$\Pi_i(\lambda_i)$ equals zero there, and it is non-zero in the
complement. This complement is connected, because $\Sing_i^j$ has
codimension 2  (by proven part of Lemma~\ref{3} and as we
explained in Section~\ref{ma+}),  and the claim follows.

At last, let us prove that such  continuous field $v_i$ exists in
the small disk neighborhood $U^n\ni x$. Since
$U^n\setminus\Sing_i^j$ is connected, we can define $v_i$ in one
of two possible ways at some point $x_0$ and extend by continuity
along paths in $U^n\setminus\Sing_i^j$.   We need to show that the result is
well-defined.

In order to do  this we connect two paths $\phi_0,\phi_1$ from
$x_0$ to $x_1$ in $U^n\setminus\Sing_i^j$ by a homotopy
$\phi_\tau$ in $U^n$. The paths and the homotopy can be assumed
smooth.  Since $\Sing_i^j$  has codimension 2,  we can perturb
homotopy and make it to  be transversal to $\Sing_i^j$. Thus, the
intersection of $\textrm{Image}_{\phi_\tau}$ with $\Sing_i^j$ is a
finite set $\{(t_k,\tau_k)\}\in[0,1]\times [0,1]$ and it suffices
to consider only one point of intersection
$y_0=\phi_{\tau_0}(t_0)=\phi(t_0,\tau_0)\in\Sing_i^j$. If we can
find the required field $v_i$ on a  transversal 2-dimensional disk
at $y_0$, we are done.

As we explained in  Section~\ref{ma+}, at  almost every point
$y\in \Sing_i^j$ we have  $\l_{i-1}(y)\ne\l_{i+1}(y)$. (Actually,
for $j=3$ this is true at every point.) Thus, without loss of
generality, we can assume that $\l_{i-1}(y_0)\ne\l_{i+1}(y_0)$.

Assume  $\l_{i-1}(y_0)\ne\l_i=\l_{i+1}(y_0)$.  The case
$\l_{i-1}(y_0)=\l_i\ne \l_{i+1}(y_0)$ is completely analogous.

Let $A=\{i,i+1\}$. Then $y_0\in\Reg_A(M)$.  Consider the leaf
$M_A(y_0)$. This is a 2-dimensional manifold transverse to
$\Sing_i^j$ at $y_0$. The homotopy can be perturbed to have the
image locally coinciding with $M_A(y_0)$. Since $v_i\in D_A$, the
problem, thanks to Lemma \ref{M_A}, is reduced to a local
2-dimensional question on $M_A(y_0)$.

Consider the restriction of the metrics to $M_A(y_0)$. Denote by
 $L_A$  the tensor $(\ref{l})$ constructed for the restrictions of the
 metrics. We denote by $\l_A\le\l_A'$ its eigenvalues.  By Lemma~\ref{M_A}, $\l_A$ is constant, $\l_A'$  is not.
 If there exists  a (continuous) vector field $v_A$ on $M_A$ such
 that it vanishes precisely at $y_0$,  such that it is
 eigenvector of $L_A$ with eigenvalue $\l_A$, and such that its length is $\sqrt{\lambda_A'-\lambda_A}$,
  we are done. Indeed, by Lemma~\ref{M_A} the vector field $v_i$  given by
 $$\sqrt{C^{-1/3}\left|\prod_{\alpha\ne i,i+1}(\lambda_i-\lambda_\alpha)\right|}\ \ \ \  v_A, $$
where $C$ is given by (\ref{c}), satisfies the conditions
(\ref{v}).  Since $$\sqrt{C^{-1/3}\left|\prod_{\alpha\ne
i,i+1}(\lambda_i-\lambda_\alpha)\right|}$$ is  a smooth positive
function, the existence of $v_A$ implies the existence of $v_i$.

    Let us prove the existence of  such vector field $v_A$.
   At every $y\in M_A(y_0)$, $y\ne y_0$,  denote by   $l_A$  the eigenspace of $L_A$ corresponding to $\l_A$.
    Let us show that
that for every geodesic $\gamma$ on $M_A(y_0)$ passing through
$y_0$ the velocity vector  $\dot\gamma(t)$  is orthogonal  (in the
restriction of $g$) to $l_A$ at every $\gamma(t)\ne y_0$.  Indeed,
let $I_t^A$ be the one-parametric family of the integrals from
Theorem~\ref{integrability} constructed for the restrictions of
$g$ and $\bar g$ to $M_A(y_0)$. Consider the integral
$I_{\l_A}^A$. At the tangent plain to every point $z$ consider the
coordinates such that the restriction of $g$ to $M_A(y_0)$ is
given by $\diag(1,1)$ and $L_A$ is $\diag(\l_A,\l_A')$. In this
coordinates, the integral $I_t^A$ equals $(\lambda'_A-t)\xi_1^2+
(\lambda_A-t)\xi_2^2$, so that $I_{\l_A}^A$ is equal to $
(\lambda'_A-\lambda_A)\xi_1^2. $
 We see that the integral vanishes on every geodesic $\gamma$ passing
through $y_0$. Because $\l_A'(z)\ne \l_A(z)$ for $z\ne y_0$, we
 obtain that the component  $\xi_1$ of the velocity vector of 
$\gamma$ at $z$ vanishes, which means that the eigenvalue of $L_A$
corresponding to $\lambda_A$ is orthogonal to $\gamma$.

\begin{center} 
\begin{figure}[h!] \label{fig}
{{\psfig{figure=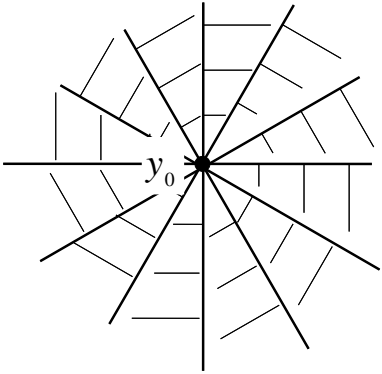}}}  
  \caption{In dimension 2, there exists a 
vector  field  orthogonal to all geodesics containing $y_0$.}\label{picture}
\end{figure}
\end{center}

Clearly,  in $M_A(y_0)\setminus y_0$ there exists a vector field
of length 1 such that  it is orthogonal to the geodesics passing
through $y_0$, see Figure~\ref{picture}.

 Multiplying  this vector field by
$\sqrt{\lambda_A'-\lambda_A}$, we obtain a required  vector field
$v_A$ on $M_A(y_0)\setminus y_0$. We put $v_A=0$ at point  $y_0$.
Since $\sqrt{\lambda_A'-\lambda_A}$ converges to 0 when $x$ tends
to $y_0$, the result is a required continuous vector field $v_A$
on $M_A(y_0)$. Therefore, there exists a vector field $v_i$ along
$M_A(y_0)$ (satisfying (\ref{v})). Thus, the vector $v_i$ at $x_1$
does not depend on the choice of path connecting $x_0$ and $x_1$.
Finally, $v_i$ is well-defined at the whole
$U^n\setminus\Sing_i^j$, and is at least continuous on it.

At the points  of  $U^n\cap\Sing_i^j$  let us put $v_i$ equal to
zero. Since $\Pi_i(\lambda_i)$ tends to $0$ when $x$ approaches
$\Sing_i^j$, the vector field is continuous on $U^n$. Statement
{\bf (S)} is proven.

 Then, at least on the double cover of $M^n$, there exists a
 continuous vector field $v_i$ satisfying (\ref{v}).
 Without loss of generality, we can assume that
the vector field $v_i$ is defined already on $M^n$.

Now let us prove that the vector  field $v_i$ is actually smooth.
Clearly, it is smooth on the compliment to $\Sing_i^j$, because it
coincides with the appropriate field $\frac{\partial }{\partial
x_i}$ there. Denote by $F_t$ the flow of the vector field $v_i$ on
$M^n\setminus(\Sing_i^2\cup\Sing_i^3)$. This flow is globally (=for
every value of $t$) defined.  Indeed, if
$x\notin\Sing_i^2\cup\Sing_i^3$, then $\lambda_{i-1}(x)<\lambda_i<
\lambda_{i+1}(x)$. Since $v_i$ is an eigenvector of $L$ with
eigenvalue $\lambda_i$ and the Nijenhuis tensor $N_L$ vanishes
(Corollary \ref{nijenhuis}), for every $t$ we have:
$\lambda_{i-1}(F_t(x))= \lambda_{i-1}(x)$, $\lambda_{i+1}(F_t(x))=
\lambda_{i+1}(x)$. Therefore,  the trajectory of the flow passing
through $x$ never approaches the set $\Sing_i^2\cup\Sing_i^3$.

The function $J(\eta):=g(v_i, \eta)$ is a linear in velocities
integral of the geodesic flow, which implies that $F_t$ acts by
isometries on $M^n\setminus(\Sing_i^2\cup\Sing_i^3)$. Since $M^n\setminus(\Sing_i^2\cup\Sing_i^3)$ is everywhere dense in $M^n$, the map $F_t$ can be
extended by completeness to act by isometries on the whole $M^n$.
Thus,  there exists a Killing vector field on $M^n$ coinciding
with $v_i$ almost everywhere. Since every Killing vector field is
smooth, the vector field $v_i$ is
smooth. Theorem~\ref{fake} is proven.

\section{Proof of Theorem~\ref{th1}}


\label{last1}

We use induction by the dimension. If dimension of the manifold is
$n<2$, Theorem~\ref{th1} is trivial. Assume that for every
dimension less than $n$ Theorem~\ref{th1} is true and consider
$\dim M=n$.

Vanishing of the topological entropy for the lift of a dynamical
system to a finite cover (of a closed manifold) implies vanishing
of the topological entropy of the original system. Thus, we assume
that already on $M^n$ for  every constant eigenvalue $\lambda_i$
we can associate a global vector field $v_i$ from
Theorem~\ref{fake}. Therefore for every constant $\l_i$ we
globally define the integral $J_i$ such that
its  differential does not vanish over the  points of $\Reg(M^n)$, 
it commutes  with all  integrals $I_t$, 
it is  functionally   dependent with the integral $I_{\lambda_i}$.

By geodesic flow we will understand the restriction of the Hamiltonian system 
on $TM^n$  with the Hamiltonian $H(\xi):=g(\xi,\xi)$ to 
$T_1M^n=\{\xi\in TM^n: \  H(\xi) =1\}$. The symplectic form on $TM^n$ came from $T^*M^n$ via standard identification by $g$.

Since  $T_1M^n$ is compact, 
 the variational principle (see, for example, Theorem 4.5.3 of \cite{KH})
 holds, and we obtain 
 $$
\h(g)=\sup_{\mu\in \mathfrak{B}}h_\mu(g).
 $$
Here $\mathfrak{B}$ is the set of all invariant ergodic
probability measures on $T_1M^n$  and $h_\mu$ is
the entropy of an invariant measure $\mu$. Recall that a measure
is called {\bf ergodic}, if $\mu(B)(1-\mu(B))=0$ for all
$\mu$-measurable invariant Borel sets $B$.

Therefore, in order to prove Theorem~\ref{th1}, it is sufficient 
 to prove that $h_\mu(g)=0$ for all $\mu\in \mathfrak{B}$.
Fix one such measure and let $\Supp(\mu)$ be its support (the set
of $x\in M^n$ such that every  neighborhood
$U_\epsilon(x)$ has  positive measure). 

Since the measure is ergodic, its support lies on a level surface
of every invariant continuous function. Then,
$\op{Supp}(\mu)$ is included into a Liouville leaf $\Upsilon$
(Recall that {\it a Liouville leaf}  is  a connected component of
the set $\{I_{t_1}=c_1,\dots,I_{t_n}=c_n\}$, where $c_1,...,c_n$
are  constants.)

 Suppose a point $\xi\in \op{Supp}(\mu)$ is
 nonsingular, or
is a  removable singular point (in the sense that every
$I_{\lambda_i}$ such that $dI_{\lambda_i}=0$ can be  replaced  by
a linear integral $J_i$ such that   $dJ_i\ne 0$). Then, a small
neighborhood $U(\xi)$ of $\xi$ in $\op{Supp}(\mu)$
\begin{itemize}
\item has positive measure in $\mu$,
\item contains only points that are  nonsingular or
removable-singular.
\end{itemize}
We will show that these two conditions  imply that the entropy of
$\mu$ is zero.

By implicit function Theorem,  $\Upsilon$ is $n$-dimensional near
$\xi$. Denote  by $O(\xi)$ the orbit of the Poisson action of $(\R^n,+)$
containing $\xi$. Since it is also $n$-dimensional, in a small
neighborhood of $\xi$ it coincides with $\Upsilon$. Thus,
$U(\xi)\subset O(\xi)$.

The orbits of the Poisson action and the dynamic on them are
well-studied  (see, for example, \S 49 of \cite{A}).  There exists
a diffeomorphism to $$ T^k\times
\mathbb{R}^{n-k}=\underbrace{S^1\times ...\times S^1}_{k}\times
\underbrace{ \mathbb{R}\times...\times \mathbb{R}}_{n-k}$$ with
the standard coordinates $\phi_1,...,\phi_k\in (\mathbb{R}\ \
\textrm{mod} \ \ 2\pi)$, $t_{k+1},...,t_n\in \mathbb{R} $ such
that in these coordinates (the push-forward of) every trajectory
of the geodesic flow is given by the formula  $$
(\phi_1(\tau),...,\phi_k(\tau),t_{k+1}(\tau),...,t_n(\tau))=(\phi_1(0)+\omega_1\tau,...,\phi_k(0)+\omega_k\tau,
t_{k+1}(0)+\omega_{k+1}\tau,...,t_{n}(0)+\omega_{n}\tau),$$ where
the constants $\omega_1,...,\omega_n$ are universal on  $T^k\times
\mathbb{R}^{n-k}$.

We see that if at least one of the constants
$\omega_{k+1},...,\omega_n$ is not zero, every point  of $U(\xi)$  
 is wandering  in  $\Supp(\mu)$ 
 (see \S3 in  Chapter 3 of    \cite{KH} for definition), 
which contradicts the invariance of the measure.  
Then, the entropy of $\mu$ is zero.

If all constants $\omega_{k+1},...,\omega_n$ are zero, the coordinates $t_{k+1},...,t_{n}$ are constants on the trajectories of the geodesic flow. Since $\mu$ is ergodic, they are constant on the points of $\Supp(\mu)$. Then, $\Supp(\mu)$ 
is (diffeomorphic to) the torus $T^{\bar k}$ of dimension $\bar k\le k$, and 
the dynamics on  $\Supp(\mu)$  is (conjugate to)  the linear  flow  on  $T^{\bar k}$.  Then, the
entropy of $\mu$ is zero, see for example Proposition~3.2.1 of \cite{KH}.

Now suppose that $\Supp(\mu)$ contains only singular points which
are not removable. If all of them belong to
$\cup_{i,j}T\Sing_i^j$, then (because the measure is ergodic)
$\Supp(\mu)$ is a subset of a certain $T\Sing_i^j$. Since
$\Sing_i^j$ is totally geodesic, and since by induction hypothesis
the topological entropy on $\Sing_i^j$ is zero, the entropy of
$\mu$ is also zero.

The last case is when $\Supp(\mu)$ contains a   singular point
which is  not removable and which does not belong to
$\cup_{i,j}T\Sing_i^j$.    Then, since all $\Sing_i^j$  are
totally geodesic, and since there are finitely many of them,
$\Supp(\mu)$ contains a   singular point $\xi$ which is  not
removable and such that its projection does not belong to
$\cup_{i,j}\Sing_i^j$. Then, the projection of a small
neighborhood $U(\xi)\subset\Supp(\mu)$ 
of $\xi$ does not contain points of
$\cup_{i,j}\Sing_i^j$.   

  From Theorems~\ref{singular},\ref{fake} it follows,
that for certain $\bar \lambda_i$ such that $\lambda_i$ is not constant 
  the differentials of  $I_{\bar \lambda_i}$ 
 vanish  at $\xi$. Since the number of such $\bar \lambda_i$ is finite, and  
 since the measure is ergodic, we obtain that there exists $i$ such that
\begin{itemize} 
\item $dI_{\bar \lambda_i}=0$ at every point of $\Supp(\mu)$, 
\item  the eigenvalue $\lambda_i$ satisfies the assumptions of
Lemma~\ref{compact}. (Otherwise the singularity is removable or 
$\xi$ lies in $\cup_{i,j}T\Sing_i^j$.)
\end{itemize} 
 Hence, by Lemma~\ref{compact},  for every
point $y$ from the projection of $U(\xi)$ we have that 
$M_{\textsf{C}(i)}(y)$ is compact.  Then, by
Theorem~\ref{singular}, for every $\eta\in U(\xi)$, 
the projection of  the trajectory
of the geodesic flow passing through $\eta$  stays on the
corresponding $M_{\textsf{C}(i)}$. Since all $M_{\textsf{C}(i)}$
passing through the projection of $U(\xi)$ 
are compact and do not intersect one another,  a  trajectory
staying in  one $T_1M_{\textsf{C}(i)}$ never approaches another
$T_1M_{\textsf{C}(i)}$. Thus, since $\mu$ is ergodic, all points of
$\Supp(\mu)$ belong to a  certain $T_1M_{\textsf{C}(i)}(y)$. Then,
the dynamics on $\Supp(\mu)$ is a subsystem  of the geodesic flow
for the  restriction of $g$ to $M_{\textsf{C}(i)}(y)$. 
(Indeed, if a geodesic of 
a metric lies on a submanifold, then it is a geodesic in the 
restriction of the metric to the submanifold.) 
Finally, by induction assumptions, 
the entropy of $\mu$ is zero.

Thus, for every ergodic probabilistic invariant measure $\mu$ its
entropy is zero. Finally, the topological entropy is zero.
Theorem~\ref{th1} is proven.


\section{Topological restrictions for manifolds with infinite
fundamental group: announcement} \label{lastsec} 

 \begin{theorem}\po \label{th2}
Suppose the  Riemannian metrics $g$ and $\bg$ on a closed
connected manifold  $M^n$ are  geodesically equivalent and
strictly non-proportional at least at one point. Then  some finite
cover of  $M^n$ is diffeomorphic to  the product $Q^k\times
T^{n-k}$ of a rational-elliptic manifold and the torus.
 \end{theorem}
The proof of this theorem  is lengthy and will appear  elsewhere (
for small dimensions, in view of Theorem~\ref{th1},  
Theorem~\ref{th2}  follows from \cite{PP2}).
Here we sketch the proof only. It  uses  Corollary~\ref{crl1},
methods developed in \cite{M1,M4} and classical results of
\cite{CG}.

 In
\cite{M1}, it was shown that if a manifold  with non-proportional
geodesically equivalent metrics has an infinite fundamental group,
it admits a local product structure (= a new Riemannian metric and
two orthogonal foliations  of complementary dimensions $B_k$ and
$B_{n-k}$ such that in a small neighborhood of  almost every point
all three object look as they come from the Riemannian product of
two Riemannian manifolds). In \cite{M4} (see Lemma~2 there), it
was shown that (assuming that the initial metrics $g$ and $\bar g$
are strictly non-proportional at least at one point), the
restriction of the local-product metric to the leaves of the
foliations admits a metric which is geodesically equivalent to it
and strictly non-proportional to it at almost every point. By
applying the same construction to the leaves, we obtain that $M^n$
admits a Riemannian metric $h$ and $m$ orthogonal foliations
$B_{k_1},B_{k_2},...,B_{k_m}$  of complementary dimension
$k_1+k_2+...+k_m=n$  such that
\begin{itemize}
\item the restriction of the metric $h$  to $B_{k_1}$ is  flat,
\item the leaves of $B_{k_2},B_{k_3},...,B_{k_m}$ are compact and have finite fundamental group
 (this is actually the lengthy part of
 the proof; its proof it similar to the proof
of Theorem~2 from \cite{M1}, but one can not apply Theorem~2 from
\cite{M1}  directly and should essentially repeat all steps  of
its proof in a slightly different setting.)
\item the restriction of $h$ to each of  $B_{k_2},B_{k_3},...,B_{k_m}$ admits
a metric which is geodesically equivalent to it and is  strictly
non-proportional to it at least at one point.
\item locally,  in a neighborhood of every point, the metric $h$ and
the foliations $B_{k_i}$ look as they (simultaneously) came from
the direct  product of $m$ Riemannian manifolds.

\end{itemize}
Then, by Corollary~\ref{crl1},  the universal cover of
$B_{k_2}\times B_{k_3}\times ... \times B_{k_m}$ is rational
elliptic, and Theorem~\ref{th2} follows from Theorem~9.2 of
\cite{CG}.

\section{ Vanishing of the entropy pseudonorm: announcement} \label{last3}

                     An action
$\Phi:(\R^n,+)\to\op{Diff}(W)$ determines the following 
{\bf entropy pseudonorm} \cite{K} :
 $$
\rho_\Phi(v):=\h(\Phi(v)).
 $$
The triangle inequality is based on the Hu's formula \cite{H}.

In particular, for the Poisson action
$\Phi:(\R^n,+)\to\op{Symp}(W^{2n},\omega)$ associated with a
Liouville-integrable Hamiltonian system one gets
a certain pseudonorm $\rho_\Phi:\R^n\to\R_{\ge0}$.

This pseudonorm is degenerate for most  examples  of
integrable geodesic flows with positive entropy ($W^{2n}=TM^n$),
but it is possible to construct a Liouville-integrable
Hamiltonian system such that  $\rho_\Phi$ is  a norm \cite{K}.

 \begin{theorem}\po\label{pseudn}
Suppose the Riemannian metrics $g$ and $\bg$ on a closed connected
manifold $M^n$ are  geodesically equivalent and 
strictly non-proportional at least at one point. 
Let $\Phi$ be the Poisson action
constructed by the integrals $I_{t_1},\dots,I_{t_n}$,   where the numbers
$t_i$ are mutually different. Then,  $\rho_\Phi(v)=0$ for every $v\in R^n$.
 \end{theorem}
The proof of this theorem will be published elsewhere.


 \hspace{-20pt} {\hbox to 12cm{ \hrulefill }}
\vspace{5pt}

{\footnotesize \hspace{-10pt}\textsc{Institute of Mathematics and
Statistics, University of Troms\o, Troms\o\ 90-37, Norway}.

\hspace{-10pt} E-mail address: \quad kruglikov@math.uit.no}
\vspace{-1pt}

\

{\footnotesize \hspace{-10pt}\textsc{Mathematisches Institut der
Albert-Ludwigs-Universit\"at, Eckerstra}\ss\textsc{e-1, Freiburg
79104, Germany}.

\hspace{-10pt} E-mail address: \quad
matveev@email.mathematik.uni-freiburg.de} \vspace{-1pt}

\end{document}